%% file: main.tex
\documentclass[review,10pt]{elsarticle}
\usepackage[usenames]{color,colortbl}
\usepackage{fancyvrb} 
\usepackage{bm}
\usepackage{amsmath}
\usepackage{subfigure}

\usepackage{multirow}
\usepackage{graphicx,graphics}
\usepackage{amssymb}
\usepackage{array}
\usepackage{url}

\usepackage{rotating}
\usepackage{multicol,lipsum}

\usepackage{listings}
\definecolor{mygreen}{RGB}{28,172,0} 
\definecolor{mylilas}{RGB}{170,55,241}
\definecolor{LightCyan}{rgb}{0.88,1,1}
\definecolor{Gray}{gray}{0.9}

\begin{document}
\begin{frontmatter}
\title{A Novel Approach for Optimal Trajectory Design with Multiple Operation Modes of Propulsion System, Part 2}

\author{Ehsan Taheri \fnref{label1f}}  
\fntext[label1f]{Assistant Professor, Email: etaheri@auburn.edu}
\address{Department of Aerospace Engineering, Auburn University\\ Auburn, AL 36849}
\author{John L. Junkins \fnref{label2f}}
\fntext[label2f]{Distinguished Professor, Email: junkins@tamu.edu}
\address{Department of Aerospace Engineering, Texas A\&M University\\  College Station, TX 77843-3141}
\author{Ilya Kolmanovsky \fnref{label3f}} \author{Anouck Girard \fnref{label4f}}
\fntext[label3f]{Professor, Email: ilya@umich.edu}
\fntext[label4f]{Associate Professor, Email: anouck@umich.edu}
\address{Department of Aerospace Engineering, University of Michigan\\  Ann Arbor, MI 48109-2140}

\begin{abstract}
Equipping a spacecraft with multiple solar-powered electric engines (of the same or different types) compounds the task of optimal trajectory design due to presence of both real-valued inputs (power input to each engine in addition to the direction of thrust vector) and discrete variables (number of active engines). Each engine can be switched on/off independently and ``optimal'' operating power of each engine depends on the available solar power, which depends on the distance from the Sun. Application of the Composite Smooth Control (CSC) framework to a heliocentric fuel-optimal trajectory optimization from the Earth to the comet 67P/Churyumov-Gerasimenko is demonstrated, which presents a new approach to deal with multiple-engine problems. Operation of engine clusters with 4, 6, 10 and even 20 engines of the same type can be optimized. Moreover, engine clusters with different/mixed electric engines are considered with either 2, 3 or 4 different types of engines. Remarkably, the CSC framework allows us 1) to reduce the original multi-point boundary-value problem to a  two-point boundary-value problem (TPBVP), and 2) to solve the resulting TPBVPs using a single-shooting solution scheme and with a random initialization of the missing costates. While the approach we present is a continuous neighbor of the discontinuous extremals, we show that the discontinuous necessary conditions are satisfied in the asymptotic limit. We believe this is the first indirect method to accommodate a multi-mode control of this level of complexity with realistic engine performance curves. The results are interesting and promising for dealing with a large family of such challenging multi-mode optimal control problems.  
\end{abstract}

\end{frontmatter}

\input{Introduction.tex}

\input{ProblemFormulation.tex}

\input{P2.tex}

\input{Conclusion.tex}

\section{Acknowledgment}
Ehsan Taheri and John Junkins would like to thank our sponsors: AFOSR (Dr. Stacie Williams) and AFRL (Dr. Alok Das) for their support and collaborations under various contracts and grants. Ehsan Taheri and Ilya Kolmanovsky would also like to acknowledge the support of the National Science Foundation under the Award Number CNS 1544844. 

\bibliographystyle{elsarticle-num}
\bibliography{References.bib}

\end{document}

%% file: Introduction.tex
\section{Introduction} \label{sec:intro}
Over the past six decades, significant technological developments have been made that drastically improve, and arguably, revolutionize the capacity of space probes to accomplish missions beyond what had been possible using earlier technologies \cite{rayman2000results,racca2003new,ceriotti2011hybrid,oh2014solar,duchek2015solar,li2019geostationary}. Specifically, the recent breakthroughs that have taken place in advanced Solar Electric Propulsion (SEP)  systems \cite{lev2019technological} are noticeable in reducing launch mass and enabling inexpensive missions \cite{rayman2002design,kuninaka2011lessons,rayman2006dawn,dachwald2008main,konstantinov2017method,sarli2017hayabusa,song2019solar}. 

Electric engines (thrusters) operate at a higher level of mass efficiency compared to chemical rockets, which leads to delivering larger payloads and ability to reach a diverse set of targets, specifically, for small-body rendezvous missions. Efficiencies can be gained not only through an order of magnitude reduction of propellant required, but also through a reduction in the mass of engine, propellant tanks, and support structure. A simultaneous advancement in low size, weight and power sensors, communication and computer systems have further expanded possibilities. As a consequence, much more can be accomplished with smaller and less expensive spacecraft. SEP is now envisioned for many in-space logistics supply purposes and for cargo missions \cite{manzella2014high,jagannatha2018optimization,woolley2019cargo}. 

In spite of the progress made recently, certain challenges still remain that have to be surmounted when SEP systems are used as the primary means of producing propulsion force. Due to very small thrust values that these engines produce, they have to operate for longer duration. It is also possible to configure a cluster of thrusters to create additional thrust and reduce the flight time. A more complicated facet of designing trajectories with SEP systems is the inherent coupling between the available power and the trajectory state (primarily inverse square dependence of available solar power on the distance to the Sun).

Furthermore, while we have approximate insight on aging of solar power systems and electric propulsion systems, the long mission duration lead to uncertainties in the actual SEP system performance versus the models used to optimize the planned trajectory. These considerations, collectively, demand a paradigm shift with respect to the tools needed for solving trajectory optimization problems \cite{patel2006maximizing}. In this paper, we deal with deterministic hybrid system models for SEP trajectory optimization. 

Trajectory optimization problems, in particular, and optimal control problems, in general, are solved by direct methods, indirect methods \cite{betts1998survey, conway2012survey} or variants thereof \cite{yang2007earth,guo2012homotopic}. The resulting optimal trajectories play a pivotal role in efficient operation of a number of flight vehicles \cite{kamyar2014aircraft, davoudi2019quad}. Here, we focus on indirect optimization methods since satisfaction of the resuling first-order necessary conditions of optimality guarantees at least a local extremal solution \cite{trelat2012optimal}. Moreover, the obtained extremals are usually of higher accuracy compared to their direct approximations. Depending on the type of problem and constraints it may be necessary to solve either a two-point boundary-value problems (TPBVPs), or in more complicated cases, a multi-point boundary-value problems (MPBVP).

For initial analysis, it is a common practice to consider simplified models for SEP systems. A single engine is assumed to operate at  \textit{effective} constant specific impulse and efficiency values \cite{bertrand2002new, jiang2012practical,tang2015low,taheri2016enhanced,junkins2018exploration}. However, the actual performance of SEP systems depends on the input power, which has to be taken into consideration for obtaining more realistic trajectories \cite{sauer1997solar,nah2001fuel,casalino2004optimization,ross2005roadmap, mengali2005fuel,petukhov2012method,oh2013simple,kluever2014efficient,zhang2014fuel,laipert2015automated}. More accurate knowledge of the capability of a spacecraft to change its trajectory is obtained when more realistic models of the SEP system and dynamics are used \cite{quarta2011minimum,chi2018power,ellison2018application}. 

Depending on the mission design stage (i.e., a preliminary design study or a final study for actual guidance), a variety of tools are developed that incorporate various levels of fidelity for sub-system models, perturbation models, and planetary ephemerides data. A fairly comprehensive review of the models, objective functions, and solution approaches commonly used for spacecraft trajectory optimization is provided in \cite{shirazi2018spacecraft}. Low- to medium- to high-fidelity models/tools for designing spacecraft trajectories exist that formulate and solve the resulting optimal control problem by various optimization methodologies \cite{sauer1997solar,sims1999preliminary,whiffen2006mystic,englander2012automated,williams2010overview}. 

The focus of this work is to demonstrate the application of the Composite Smooth Control (CSC) framework developed in part 1 \cite{taheri2019anovelpart1} to problems with multiple modes of operation. In particular, we study fuel-optimal interplanetary trajectory optimization when the propulsion system consists of a cluster of engines. In the case of multiple engines, we permit multiple modes in the sense that one, or two or many engines can be selected in such a fashion that optimality conditions are satisfied. A prominent feature of the proposed approach is that not only the optimal instances of transition between different operating modes, but also the optimal number of engines as well as their operating conditions are revealed without \textit{a priori} assumptions. 

The key difference in the CSC framework (compared to Ref. \cite{saranathan2017relaxed}) is in simultaneous smooth transitions between possible modes of operations and having smooth bang-bang control inputs. The former is achieved by constraints (power-driven constraints in the considered problems), whereas the latter is achieved by the switching function of the associated control mode. In the CSC framework, bang-bang control profiles of thrusters are also incorporated. There is an important implementation subtlety when formulating and solving problems that consist of bang-bang type controls: \textit{the constraint that determines the activation of the bang-bang control is the so-called switching function associated with that specific control}. In other words, for each control input that has a bang-bang control structure, the switching function, $S$, serves as the constraint ($S = 0$) and is considered as the argument of the hyperbolic tangent smoothing (HTS) method \cite{taheriagenericsmoothing2018}. Therefore, the time associated with this constraint is obtained implicitly, and in an autonomous manner. This means that the time of control switches (between the two extreme limits) are also determined autonomously. Collectively, in addition to incorporating smooth transitions due to multiple time- or state-triggered constraints, CSC enables us to incorporate switching-function-triggered constraints that govern bang-bang type control inputs.

%% file: ProblemFormulation.tex
\section{Power System and Perturbation Modeling} \label{sec:models}
In this section, a review of the solar electric power models, and thruster performance data is given. Then, equations of motion are given while taking into account the power sub-system and actuation models.

\subsection{Solar Arrays and Spacecraft Sub-system Power Models}
There are multiple ways to model power sub-system of spacecraft \cite{sauer1997solar,quarta2011minimum}. A detailed discussion of the power sub-system is given in Part 1 and we review the main points. The available power, $P_{\text{ava}}$, to be distributed to engines can be calculated as
\begin{equation} \label{eq:avapower}
P_{\text{ava}} = P_{\text{SA}}(t,r) - P_{s/c}(r),
\end{equation}
where 
\begin{align}
P_{\text{SA}}(t,r) & = \psi(t) \phi(r) P_{0,\text{BOL}},\\
\phi(r) & = \frac{1}{r^2} \left [ \frac{A_1+\frac{A_2}{r}+\frac{A_3}{r^2}}{1+A_4 r +A_5 r^2} \right ], \label{eq:powerdist}\\
\psi(t) & = (1-\sigma)^{\tau(t)}, \label{eq:powertime}
\end{align}
where $P_{0,\text{BOL}}$ denote the nominal beginning-of-life power produced by the solar arrays at one astronomical unit (AU) from the Sun. $P_{s/c}(r)$ denotes the power needed for operation of spacecraft systems and the power processing unit (PPU). In Eq.~\eqref{eq:powertime}, $\sigma$ denotes the efficiency decay rate of the solar arrays measured as a fractional rate per year (usually 0.02 to 0.04 per year) and $\tau(t)$ denotes elapsed time from the launch time measured in years \cite{ellison2018application}. For numerical results in this paper, $P_{s/c} = P_{\text{ppu}}$ and $P_{\text{ppu}}$ is assumed to be the required power to operate all spacecraft sub-systems other than the engine(s). The empirical approximations of Eq.~\eqref{eq:powerdist} is fit to represent a more complicated model; the details and the $A$ coefficients are given in Part 1 \cite{taheri2019anovelpart1}. While a single PPU is able of powering a limited number of engines \cite{oh2014solar}, we assumed that only one PPU is capable of powering all of the engines in the numerical results. The proposed framework can readily extend to situations that multiple PPUs have to be used, but, our focus is on the utility of the general CSC framework.   

\subsection{Practical Engine Models}
\begin{figure}[htbp!]
\centering
\includegraphics[width=4.6in]{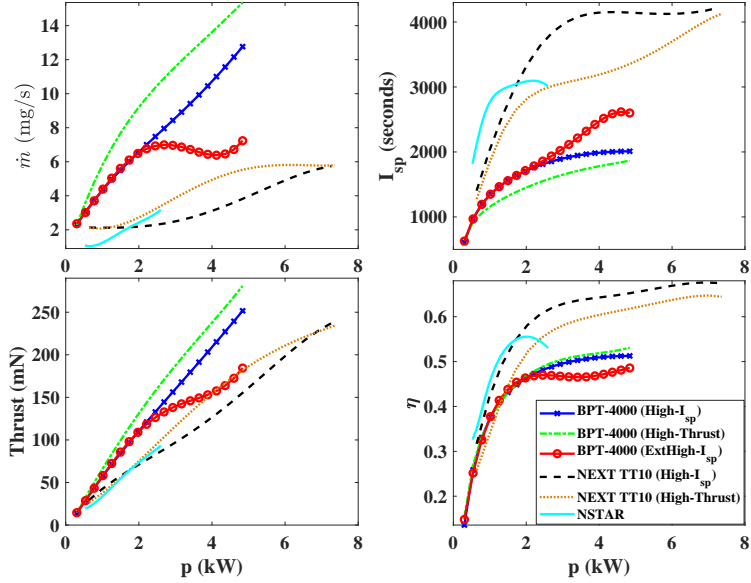}
\caption{Performance of engines listed in Table \ref{tab:enginepar} over their admissible power ranges \cite{ellison2018application,hofer2010high}.}
\label{fig:powerpercur}
\end{figure}

In low-fidelity trajectory design, the dependency of both specific impulse and thrust value to power are ignored, and instead, \textit{effective} constant values are used. In this work, however, we consider more realistic engine models. These models represent approximations of the actual performance of SEP engines. These are based on experiments and simulations that cover a finite, but a large set (on the order of 40 to 100 or larger) of operating points. Each operating point is characterized by a local mass flow rate and thrust magnitude \cite{rayman2002design}. A best-fit second-order polynomial interpolation approximation is shown to be an adequate surrogate model for a preliminary mission analysis \cite{koppel2003smart,mengali2007tradeoff}. Rather than linear interpolation using a grid, continuous interpolation can be used. We adopt fourth-order polynomials interpolation (with power as the independent variable) for thrust and mass flow rate defined as
\begin{align} 
T(P_{\text{en}}) = &~a_{\text{T}} P_{\text{en}}^4 + b_{\text{T}}  P_{\text{en}}^3+ c_{\text{T}} P_{\text{en}}^2+ d_{\text{T}}  P_{\text{en}} + e_{\text{T}} ,\label{eq:quarticT}\\
\dot{m}(P_{\text{en}}) = &~a_{\text{m}}  P_{\text{en}}^4 + b_{\text{m}} P_{\text{en}}^3+ c_{\text{m}} P_{\text{en}}^2+ d_{\text{m}} P_{\text{en}} + e_{\text{m}}, \label{eq:quarticmdot}
\end{align}
where $P_{\text{en}} \in [P_{\text{min}},P_{\text{max}}]$ denotes the engine input power and $P_{\text{min}}$ and $P_{\text{max}}$ represent the lower and upper limits, respectively, over which an engine operates. The power variation interval depends on the particular engine's operational characteristics. For these representative quartic polynomial surrogate models given in Eqs.~\eqref{eq:quarticT} and \eqref{eq:quarticmdot}, Table \ref{tab:enginepar} summarizes the surrogate model coefficients of six SEP engines \cite{hofer2010high,ellison2018application}, which we adopt for the present study. The input power is in kilowatts and the mass flow rate is in milligram/s. Figure \ref{fig:powerpercur} shows the performance curves of the six engines listed in Table \ref{tab:enginepar}. 

\begin{table}[ht!] 
\begin{center} 
		\caption{Engine performance coefficients; $P$, $T$, $\dot{m}$, and $I_{\text{sp}}$ are given in kilowatts, miliNewton, milligram/s, and seconds, respectively.}\label{tab:enginepar}
		{\scriptsize
		\begin{tabular}{p{0.7cm} c c c c c c}
        \hline
        \hline
          & BPT‐4000              & BPT‐4000    & BPT‐4000                 & Next TT10            & NEXT TT11   & NSTAR \\
                     & High-$I_{\text{sp}}$  & High-Thrust & Ext‐High-$I_{\text{sp}}$ & High-$I_{\text{sp}}$ & High-Thrust & \\
                     \hline
         ID    & 1 & 2 & 3 & 4 & 5 & 6\\
        \hline
        $a_{\text{T}}$  & -0.095437 & 0.173870  & 1.174296   & -0.19082    & 0.101855017  &5.145602\\
        $b_{\text{T}}$  & 1.637023  & -1.150940 & -10.102479 & 2.96519     & -2.04053417  &-36.720293\\
        $c_{\text{T}}$  & -9.517167 & -2.118891 & 19.422224  & -14.41789   & 11.4181412   &90.486509\\
        $d_{\text{T}}$  & 72.030104 & 77.342132 & 47.927765  & 54.05382    & 16.0989424   &-51.694393\\
        $e_{\text{T}}$  & -7.181341 & -8.597025 & -1.454064  & -1.92224e-6 & 11.9388817   &26.337459\\
        \hline 
        $a_{\text{m}}$  & -0.008432 & -0.011949 & 0.086106  & -0.004776    & 0.011021367  &0.36985\\
        $b_{\text{m}}$  & 0.148511  & 0.235144  & -0.727280 & 0.05717      & -0.207253445 &-2.5372\\
        $c_{\text{m}}$  & -0.802790 & -1.632373 & 1.328508  & -0.09956     & 1.21670237   &6.2539\\
        $d_{\text{m}}$  & 3.743362  & 6.847936  & 1.998082  & 0.03211      & -1.71102132  &-5.3568\\
        $e_{\text{m}}$  & 1.244345  & 0.352444  & 1.653105  & 2.13781      & 2.75956482   &2.5060\\
        \hline
$P_{\text{min}}$ & 0.302 & 0.302 & 0.302 & 0.638 & 0.64 & 0.525\\
$P_{\text{max}}$ & 4.839 & 4.839 & 4.839 & 7.266 & 7.36 & 2.6\\
 \hline
$T_{\text{@Pmin}}$ & 13.748  &   14.537    & 14.523    & 29.356      &  26.401  & 19.215 \\
$T_{\text{@Pmax}}$ & 251.68 &   280.97    & 184.42    & 237.16       &  234.28  & 93.37\\
$\dot{m}_{\text{@Pmin}}$ & 2.305  &  2.278     & 2.358    & 2.132    &  2.11    & 1.078\\
$\dot{m}_{\text{@Pmax}}$ & 12.765 &  15.358     & 7.235   & 2.123    &  5.786   & 3.162\\
$I_{\text{sp@Pmin}}$  & 608.03 &  650.71     & 627.95    &  1404.19  &  1275.68 & 1817.05\\
$I_{\text{sp@Pmax}}$  & 2010.55 &  1865.46     & 2599.43   & 4217.93 &  1839.50 & 2106.19\\
$\eta_{\text{@Pmin}}$  & 0.136 &   0.154    & 0.148    &  0.317    &  0.258     & 0.326\\
$\eta_{\text{@Pmax}}$  & 0.513 &   0.531    & 0.486   &   0.675   &  0.644      & 0.530\\
\hline
        \end{tabular}
		}
	\end{center}
\end{table}

\subsection{Equations of Motion}
In the heliocentric phase of flight, the spacecraft motion is predominantly governed by the gravitational attraction of the Sun. In addition, the spacecraft is equipped with a cluster of engines and its mass, $m$, changes due to the consumption of propellant. The equations of motion are expressed in terms of the modified equinoctial orbital elements (MEEs) \cite{walker1986set} due to them being found superior to other sets of coordinates/elements for describing low-thrust trajectories \cite{taheri2016enhanced,junkins2018exploration}. 

The derivation of the optimality conditions is achieved in a simpler manner if a distinction is made between the trajectory dynamics of MEEs and time rate of change of spacecraft mass. Let $\textbf{x} = [p,f,g,h,k,l]^{\top}$ denote the vector of MEEs and let $\textbf{a} = [a_r, a_t, a_n]^{\top}$ represents all of the non-two-body gravitational accelerations expressed in the LVLH frame acting on the spacecraft. The dynamics of MEEs can be written as
\begin{align} \label{eq:meeddynamics}
\dot{\textbf{x}} =& \textbf{A}(\textbf{x},t) + \mathbb{B}(\textbf{x},t) \textbf{a},
\end{align}
where $\textbf{A} \in \mathbb{R}^{6 \times 1}$ denotes the two-body gravity-induced part of the dynamics and $\mathbb{B} \in \mathbb{R}^{6 \times 3}$ is the control influence matrix (these matrices are defined in \cite{junkins2018exploration}).
The MEE set has five slow variables and one (very regular) fast variable, $l$. The acceleration vector, $\textbf{a}$, in Eq.~\eqref{eq:meeddynamics} can be written as 
\begin{equation}
\textbf{a} = \textbf{u}_{\text{prop}} + \textbf{a}_{\text{sb}},
\end{equation}
where $\textbf{u}_{\text{prop}}$ denotes the control acceleration vector due to the engines and $\textbf{a}_{\text{sb}}$ denotes the perturbing accelerations vector due to secondary bodies (i.e., the other planets in the Solar System or solar radiation pressure). Another state of the system is the spacecraft mass and its time rate of change depends on the operation of engines. The details of modeling thrust acceleration term due to operation of engines, $\textbf{u}_{\text{prop}}$, and mass time rate of change are explained in the next section. 

\section{Propulsion System with Multiple SEP Engines} \label{sec:variableIspvariablethrust}
For a cluster of engines, the complete dynamical system consists of the time rate of change of spacecraft mass. Let $N_e$ denote the number of engines, and under the assumption that engine thrust vectors are all aligned, the total thrust vector produced by the engines, $\textbf{T}$, and time rate of change of spacecraft mass, $\dot{m}$ can be written as
\begin{align} \label{eq:sumthrust}
\textbf{T} & =  \hat{\bm{\alpha}} \sum_{i = 1}^{N_e} T_i, & \dot{m} = & \sum_{i=1}^{N_e} \dot{m}_i,
\end{align}
where $\hat{\bm{\alpha}}$ denotes the unit thrust steering vector and $T_i$ and $m_i$ denote the thrust magnitude and time rate of change of mass due to the operation of the $i$-th engine (for $i = 1,\cdots, N_e$). For simplicity, we assume all thrust vectors are parallel (direction $\hat{\bm{\alpha}}$). In practice, cant angles are frequently used such that all thrust vectors act along a line that passes through the mass center.  In the following, we investigate two types of clusters: 1) an engine cluster with the same type of engines, and 2) an engine cluster with different (or mixed) types of engines. The first cluster type is considered initially to explain the formulation and structure of the problem. The second type of cluster is a bit more involved, but it can be viewed as an extension of the first type of cluster.

\section{Formulation of Fuel-Optimal Boundary-Value Problem} \label{sec:formualtionofOPC}
For a fuel-optimal problem, the cost functional is written in Mayer form as 
\begin{eqnarray}
J = -m(t_f),
\end{eqnarray}
where the minimization of the negative value of the final mass is sought. For an initial fixed dry mass, the above objective corresponds to minimizing the propellant consumption. Let $\bm{\lambda} = [\lambda_p,\lambda_f,\lambda_g,\lambda_h,\lambda_k,\lambda_l]^{\top}$ denote the costate vector associated with the MEEs, $\textbf{x}$ and let $\lambda_m$ denote the costate associated with the mass state, $m$. The Hamiltonian becomes
\begin{equation} \label{eq:VIVTHamil}
H = \bm{\lambda}^{\top} \left [ \textbf{A}(\textbf{x},t) +  \mathbb{B}(\textbf{x},t) \left [ \frac{\textbf{T}}{m} + \textbf{a}_{\text{sb}} \right ] \right ] - \lambda_m \dot{m}.
\end{equation}

From the expression for the Hamiltonian given in Eq.~\eqref{eq:VIVTHamil}, together with the Euler-Lagrange equation, the dynamics of the costates can be derived as
\begin{align}
\dot{\bm{\lambda}} = & -\left [ \frac{\partial H}{\partial \textbf{x}} \right ]^{\top}  \nonumber, & \dot{\lambda}_m = & -\frac{\partial H}{\partial m}.
\end{align}

We should mention that in all our recent and current research, we make frequent use of an automated code (that is developed in MATLAB and employ the symbolic toolbox features) to derive the costate dynamics associated to the state dynamics. For instance, in this work, the contribution of the two-body gravitational model to costate dynamics,  $\dot{\bm{\lambda}}_{\text{two-body}}$, is calculated through the automated tool (or derived algebraically by hand). However, the contribution of the secondary body gravity perturbations, $\dot{\bm{\lambda}}_{\text{sb}}$, are evaluated through the complex-based derivative method, and are added numerically to the costate dynamics associated with the two-body dynamics (see Part 1 \cite{taheri2019anovelpart1} for details). Consequently, the total contribution due to secondary bodies can be written as
\begin{equation}
\dot{\bm{\lambda}}_{\text{sb}} = \sum_{i = 1}^{N_{\textbf{sb}}} \dot{\bm{\lambda}}_{\text{sb},i},
\end{equation}
where $N_{\textbf{sb}}$ denotes the number of considered secondary bodies. The costate dynamics associated with the MEEs can then be computed as
\begin{equation}
\dot{\bm{\lambda}} = \dot{\bm{\lambda}}_{\text{two-body}} + \dot{\bm{\lambda}}_{\text{sb}}.
\end{equation}

Depending on the type of maneuver, different boundary conditions can be enforced. In this paper, we are dealing with time-fixed rendezvous-type problems, where only the final mass is free. The final position and velocity vectors of the spacecraft in these types of maneuvers are required to match their target body counterparts (denoted by subscript `\textit{T}'). It is assumed that the spacecraft leaves the Earth's Sphere of Influence (SOI) on a parabolic trajectory, i.e., with zero hyperbolic excess velocity $v_{\infty} = 0$ and that, on a solar scale, the SOI is negligibly small compared to 1 AU. This is a frequently used first approximation in the solar orbit transfer preliminary mission design. The final boundary conditions (seven equality constraints) can be written in the vector function form as
\begin{equation}\label{eq:posvelconp1}
\bm{\psi}(\textbf{x}(t_f),t_f) = \left [ [\textbf{x}(t_f) - \textbf{x}_T]^{\top},\lambda_m(t_f) + 1 \right ]^{\top} = \textbf{0},
\end{equation}
where $\textbf{x}_T = [p_T,f_T,g_T,h_T,k_T,l_T]^{\top}$ denotes the target state values. Since a rendezvous maneuver is considered, the final value of the mass costate has to be -1 (due to Mayer problem transversality condition). Another unknown of the problem, for multi-revolution trajectories, is the number ($N_{\text{rev}}$) of en-route revolutions  in the transfer orbit, which has to be determined. Its value is taken into consideration when the change in the true longitude, $l$, is to be enforced as a boundary condition,
\begin{equation}
l_T = l_f + 2 \pi N_{\text{rev}},
\end{equation}
where $l_f$ is the true longitude of the final point (corresponding to $N_{\text{rev}} = 0$). The final value, $l_T$, is the updated target value for the final true longitude taking into account any additional number of revolutions. Previous studies \cite{taheri2018how,taheri2018unified} have shown that for each number of revolutions in the feasible set of $N_{\text{rev}}$ integers, there is one local extremal for fuel-optimal trajectories.  

Let $\textbf{z} = [\textbf{x}^{\top},m,\bm{\lambda}^{\top},\lambda_m]^{\top}$ denote the state-costate vector, then, we can write,
\begin{equation} \label{eq:F}
\dot{\textbf{z}} = \textbf{F} = [\dot{\textbf{x}}^{\top}, \dot{m}, \dot{\bm{\lambda}}^{\top}_{\text{two-body}}, \dot{\lambda}_m]^{\top},
\end{equation}
where $\bm{\alpha} = \bm{\alpha}^*$, $P = P^*(S)$, and $c = c^*(S)$. Once the optimal values of the control components are substituted into $\textbf{F}$, the equations of motion can be integrated numerically, if initial conditions are fully specified. 

In the type of maneuver we consider, only the initial states, $\textbf{x}(t_0) = \textbf{x}_0$ and $m(t_0) = m_0$, are specified. The final state $\textbf{x}(t_f)$ as well as the final costates are functions of the initial costates, $\bm{\eta}(t_0)$, where $\bm{\eta}(t_0) = [\bm{\lambda}^{\top}(t_0),\lambda_m(t_0)]^{\top}$ is the vector of seven unknown variables to be determined such that Eq.~\eqref{eq:posvelconp1} is satisfied. Thus, we have a TPBVP that requires a starting estimate $\bm{\eta}(t_0)$ within the domain of convergence of the particular algorithm to satisfy the prescribed boundary conditions. There are seven constraints in Eq.~\eqref{eq:posvelconp1} and seven unknown elements in \cite{taheri2016enhanced,junkins2018exploration}.

%% file: P2.tex
\section{Operation Logic for Spacecraft with Multiple SEP Engines} \label{sec:operationlogic}
This problem is quite challenging since not only the number of operating engines, but also the operating power level of active engines are not known \textit{a priori}, and have to be optimized. The traditional approach to handle these cases is to evaluate the Hamiltonian for all possible combinations of control choices and to select control inputs such that the Hamiltonian is minimized (or maximized, depending upon the formulation) \cite{jones1984advanced,quarta2011minimum}. 

There are multiple strategies to specify operational logic for multi-engine spacecraft \cite{ellison2018application}, (e.g., 1) the maximum number of active engines, and 2) the minimum number of active engines). We focus our attention to a particular case where each engine needs to be switched off or operate \textit{only} either at its maximum or minimum power setting, i.e., $P_{\text{en}} \in \{0, P_{\text{min}}, P_{\text{max}} \}$. In general, we allow for two power settings, whereas the zero-power setting is a consequence of taking into account the optimality criteria. In other words, each of the considered power settings has a corresponding switching function, which determines whether that particular power setting is active or not. The proposed methodology, as will be explained, takes into account all combinations and the optimality principle is used to select the ``best'' combination of engines per considered selection strategy.  

Many ``if $\rightarrow$ then'' conditions will emerge from the optimality conditions and trial solutions of the resulting necessary conditions have to impose these switches as the trajectory is propagated (so that the proper combination of engines and their associated power inputs are selected such that the Hamiltonian is minimized).  Consequently, we look to impose a specific set of operating conditions for the several existing engines. 

In this work, we propose a simple innovative workaround to the above problem, which involves two steps: 1) first, a classification step is performed to determine all of the possible combination of engines' operation modes in advance and out of the optimization process. This leads to a set of finite discrete operation modes (each with its own power level). 2) It is then straightforward to apply the CSC framework outlined in part 1 in order to ensure that certain operation modes are used for distribution of power through activation functions. 

An interesting aspect of using the proposed strategy is that a \textit{trivial pre-processing} calculation can be used to produce the operation modes once the number of engines and their power range of operation are determined. In practice, the number of engines would rarely exceed eight, i.e., $N_e<8$ \cite{sauer1997solar}. In addition, it is anticipated that not all distinct combinations of the operation modes will necessarily be utilized during a particular optimal maneuver. 

\section{Same-Type Engine Clusters} \label{sec:sametype}
Under the assumption that all engines are of the same type, we can express the thrust and time rate of change of mass as
\begin{align}
    \textbf{T} & =  \hat{\bm{\alpha}} \left [ T_{\text{Ne@Pmax}} + T_{\text{Ne@Pmin}} \right ],\\
    \dot{m}  & = \dot{m}_{\text{Ne@Pmax}} + \dot{m}_{\text{Ne@Pmin}},
\end{align}
where $T_{\text{Ne@Pmax}}$ ($T_{\text{Ne@Pmin}}$) denotes the thrust magnitude corresponding to the number of engines that operate at their maximum (minimum) power levels. Similarly, $\dot{m}_{\text{Ne@Pmax}}$ ($\dot{m}_{\text{Ne@Pmin}}$) denotes the time rate of change of mass corresponding to the number of engines that operate at their maximum (minimum) power levels. 

Let $N_{\text{e@Pmax}}$ and $N_{\text{e@Pmin}}$ denote the number of engines that operate at their maximum and minimum admissible power limits, respectively. We can write
\begin{align*}
 T_{\text{Ne@Pmax}} & =  N_{\text{e@Pmax}} \frac{2 P_{\text{max}}  \eta_{\text{Pmax}}}{c_{\text{Pmax}}}, & T_{\text{Ne@Pmin}} & =  N_{\text{e@Pmin}} \frac{2  P_{\text{min}} \eta_{\text{Pmin}}}{c_{\text{Pmin}}},\\
 \dot{m}_{\text{Ne@Pmax}} & =  -\frac{T_{\text{Ne@Pmax}}}{c_{\text{Pmax}}}, & \dot{m}_{\text{Ne@Pmin}} & =  -\frac{T_{\text{Ne@Pmin}}}{c_{\text{Pmin}}}.
\end{align*}

Here, $\eta_{\text{Pmax}}$ and $\eta_{\text{Pmin}}$ correspond to the efficiency of the engines at maximum and minimum power settings, respectively. Also, $c_{\text{Pmax}}$ and $c_{\text{Pmin}}$ are the maximum and minimum exhaust velocities. By substituting the thrust acceleration, $\textbf{u}_{\text{prop}} = \textbf{T}/m$ into the RHS of Eq.~\eqref{eq:meeddynamics} we have 
\begin{align} \label{eq:P2-statedynamics}
\dot{\textbf{x}} & =  \textbf{A} + \frac{1}{m} \left [ T_{\text{Ne@Pmax}} + T_{\text{Ne@Pmin}} \right ] \mathbb{B} \hat{\bm{\alpha}} , \\
\dot{m} & =    \dot{m}_{\text{Ne@Pmax}} + \dot{m}_{\text{Ne@Pmin}}.
\end{align}

Recall that each engine can only operate at two modes. Thus, the unknown control inputs that have to be determined are $\hat{\bm{\alpha}}$, $N_{\text{e@Pmax}}$, and $N_{\text{e@Pmin}}$, i.e., the direction of the total thrust vector and the number of engines operating at their respective upper and lower power limits. 

\subsection{Operation Modes For Same-Type Engine Clusters} \label{sec:impdetails}
An example is considered for explaining the process. Specifically, a low-thrust trajectory optimization from the Earth to the comet 67P/Churyumov-Gerasimenko is considered while several clusters of engines are available with BPT-4000 Extra-High-$I_{\text{sp}}$ engines. In this paper, for demonstration purposes, clusters with 4, 6, 10 and even 20 engines are considered. 

For a BPT-4000 Extra-High-$I_{\text{sp}}$ engine, $P_{\text{en}} \in [302,4839]$ Watts (see Table \ref{tab:enginepar}). It is assumed that the nominal power at 1 AU is $P_{@1AU} = 30$ KWatts and a minimum of $P_{\text{s/c}} = 500$ Watts is devoted to subsystems of the spacecraft. In order to apply the engine selection logic, we need to define the lower and upper bounds on the distance to determine the drop in nominal power. The range of distance from the Sun can be determined through analysis of the orbital elements of the target body and a minimum distance to the Sun that is specific to each problem. For this representative example, we have chosen $r \in [0.8, 2.0]$ AU. 

Table \ref{tab:representativepowerlevels} summarizes all operation modes. Here, $P_{\text{used}}$ denotes the aggregate power used by the engines for the specified mode and number of engines operating at the maximum or minimum power. We emphasize that if an operation mode is engaged, it is still possible for some of the engines in that particular operation mode not to operate. In fact, the operation of an engine is engaged by the optimality principle through the respective switching function of that particular engine. This is identical to a scenario in which there is enough power to switch an engine on, however, the optimality principle is not satisfied, which means that an engine will not switch on. 

In Table \ref{tab:representativepowerlevels}, $N_{\text{e@Pmax}}$ ($N_{\text{e@Pmin}}$) denotes the number of engines that operate at their corresponding maximum (minimum) power settings. For the prescribed values, there are 14 operating modes when $N_e = 4$. In general, the number of operation modes depends on several factors including engines power bounds, number of total engines, nominal power at 1 AU, and the prescribed range for the distance from the Sun. Also, hardware specific operating constraints may affect the number of modes one considers in a given problem. The four-engine example considered here admits all 14 mathematical permutations. We emphasize that no optimization is performed at this stage and these mode definition results are obtained through a trivial pre-processing step. 

\begin{table}[htbp!] 
\begin{center} 
		\caption{Representative operation modes for a trajectory from the Earth to comet 67P with four BPT-4000 Extra-High-$I_{\text{sp}}$ engines ($N_e = 4$); power unit is Watts. } \label{tab:representativepowerlevels}
		{\small
		\begin{tabular}{c c c c | c c c c}
        \hline
        \hline
         Mode \#         & $P_{\text{used}} $  &  $N_{\text{e@Pmax}}$  &  $N_{\text{e@Pmin}}$ &  Mode \#         & $P_{\text{used}} $  &  $N_{\text{e@Pmax}}$  &  $N_{\text{e@Pmin}}$ \\
         \hline          
          1 & 19,356      & 4      & 0    & 8  & 5443 &	1 &	2  \\
          2 & 14,819      & 3      & 1    & 9  & 5141 &	1 &	1   \\
          3 & 14,517      & 3      & 0    & 10 & 4839 &	1 &	0    \\
          4 & 10,282      & 2      & 2    & 11 & 1208 &	0 &	4 \\
          5 & 9,980       & 2      & 1    & 12 & 906  &	0 &	3   \\
          6 & 9,678       & 2      & 0    & 13 & 604  &	0 &	2     \\
          7 & 5,745       & 1      & 3    & 14 & 302  &	0 &	1   \\
        \hline
        \hline
        \end{tabular}
		}
	\end{center}
\end{table}

Note also that efficiency degradation of the solar arrays is not taken into account when we generated the data in Table \ref{tab:representativepowerlevels}. However, this will not impact the entirety of the problem. Each operation mode is characterized by its unique power value, $P_{\text{used}}$. On the other hand, efficiency degradation of the solar arrays will eventually impact the available power, $P_{\text{ava}}$, through Eq.~\eqref{eq:avapower}. At any time instant, the effects due to efficiency degradation or any other type of losses will only impact \textit{the time of transition between operation modes}, which is automatically determined during the optimization process. In the case that multiple PPUs have to operate, we just need to add the required power to $P_{\text{used}}$ and this allows us to handle the cases with multiple running PPUs. Application of the CSC framework to smoothly make a transition between operation modes is explained in the following section.

\subsection{Implementation Details for Same-Type Engine Clusters}
The most crucial step in the proposed scheme is to determine the number of operating engines, i.e., $N_{\text{e@Pmin}}$, $N_{\text{e@Pmax}}$, which are discrete values. Once these integers are known, we can evaluate the RHS of the set of state-costate dynamics and propagate them numerically to solve the shooting problem. However, we managed to overcome this difficulty by re-casting the problem into a set of finite operation modes with their associated feasible number of engines. As a consequence, if we manage to find an approach to determine the best operation mode, we automatically have information about the \textit{potential} number of engines for that operation mode. However, this does not mean that all engines within an operation mode are used; optimality criteria govern the engine operation within each operation mode.

It is at this stage that CSC framework is employed. The following steps have to be performed at each time instant, $t \in [t_0, t_f]$:
\begin{itemize}
    \item \textbf{Step 1}: $P_{\text{ava}}$ is calculated using Eq.~\eqref{eq:avapower}.
    \item \textbf{Step 2}: depending on the value of $P_{\text{ava}}$, one of the operation modes will typically be weighed more heavily compared to the other modes (and it is not known which mode will become engaged). More generally, let $\sigma \in \mathbb{N}$ denote the total number of operation modes and also let $\textbf{w} = [w_i,\cdots, w_{\sigma}]^{\top}$ (with $w_i \in [0,1]$ for $i = 1, \cdots, \sigma$) denote the vector of activation weights corresponding to the total number of operation modes. We are interested in a smooth approximation, thus, we have to \textit{sum over weighted contributions of all existing operation modes}. 
    
    Except for the operation modes \#1 and \#14, each $w_i$ consists of two multiplicative activation functions. For example, if $P_{\text{ava}}> 19,356 \text{W}$, the first operation mode is engaged, and if $P_{\text{ava}}< 302 \text{W}$, there is not enough power at all. For instance, if $ 14517 \leq P_{\text{ava}} < 14819$, then there is also enough power to activate the third operation mode. The activation weight of the third operation mode can be written in terms of two multiplicative activation functions (consistent with the convention defined in part 1 \cite{taheri2019anovelpart1}, i.e., particular constraint followed by the number of involved constraints) as
    \begin{equation*}
        w_3 = \zeta_{w_3,1} \zeta_{w_3,2},
    \end{equation*}
    where we make use of the two inequality constraints to define the smooth activation functions
    \begin{align*}
        \zeta_{w_3,1} & = \frac{1}{2} \left [1-\text{tanh}\left ( \frac{g_{w_3,1}}{\rho_c} \right)  \right], & g_{w_3,1} = P_{\text{ava}} -14819 < 0, \\
        \zeta_{w_3,2} & = \frac{1}{2} \left [1-\text{tanh}\left(
        \frac{g_{w_3,2}}{\rho_c} \right)  \right], & g_{w_3,2} = 14517 - P_{\text{ava}} \leq 0. 
    \end{align*}
    
    In the above relations, $g_{w_3,1}$ and $g_{w_3,2}$ denote the first and the second constraints that have to be taken into account for specifying the activation weight associated with $w_3$. $\rho_c$ denotes constraint-type smoothing parameter to be used during the homotopic process. The same procedure can be followed for the other activation weight functions according to the above procedure. Operation modes \#1 and \#14 can be treated similarly. It is easy to verify that in the interior of power range, $ 14,517 \leq P_{\text{ava}} < 14,819$, that $w_3 \approx 1$, and, we will find all other analogous weights are approximately zero. Additionally, as $P_{\text{ava}}$ approaches $P_{\text{used}} = 14,517 $ watts, the transition from $w_3 \approx 0$ to $w_3 \approx 1$ is smooth and an uncontrollably sharp switch from zero to unity is controlled by selection of $\rho_c$ suitably small. For all operation modes we have used the same smoothing parameter, $\rho_c$. The details follow the general procedure outlined in Section 4.1 in Part 1. The output of this step is $\textbf{w}$, which will be used in the next step.
    \item \textbf{Step 3}: Let $\tilde{N}_{\text{e@Pmax}} \in \mathbb{N}^{\sigma}$ denote the vector associated with number of engines operating at $P_{\text{max}}$, and let $\tilde{N}_{\text{e@Pmin}} \in \mathbb{N}^{\sigma}$ denote the vector associated with number of engines operating at $P_{\text{min}}$. For example, if we write Table \ref{tab:representativepowerlevels} as a three-column table, the second column represents the $\tilde{N}_{\text{e@Pmax}}$ vector and the third column represents the $\tilde{N}_{\text{e@Pmin}}$ vector. The composite smooth representations (superscript `s') of $N_{\text{e@Pmax}}$ and $N_{\text{e@Pmin}}$ can be expressed as
    \begin{align}
        N_{\text{e@Pmax}}^s & = \textbf{w}^{\top} \tilde{N}_{\text{e@Pmax}},\\
        N_{\text{e@Pmin}}^s & = \textbf{w}^{\top} \tilde{N}_{\text{e@Pmin}}.
    \end{align}
    At this stage, we have smooth, controllably sharp, representations for the problematic discrete design variables. However, we should construct smooth representations for the power sent to engines as well.  
    \item \textbf{Step 4}: First, we form the switching functions corresponding to each engine's operating power as
    \begin{align}
        S_{@\text{Pmax}} = \frac{ \| \bm{\lambda}^{\top} \mathbb{B} \|}{m} + \frac{\lambda_m}{c_{\text{@Pmax}}}, \label{eq:sfpmax}\\
        S_{@\text{Pmin}} = \frac{ \| \bm{\lambda}^{\top} \mathbb{B} \|}{m} + \frac{\lambda_m}{c_{\text{@Pmin}}}. \label{eq:sfpmin}
    \end{align}
    We can form analogous smooth thrust activation functions for power using the following relations
    \begin{align}
        \zeta_{\text{Pmax}} &  = \frac{1}{2} \left[1+\text{tanh}\left (\frac{S_{@\text{Pmax}}}{\rho_b} \right) \right ],\\
        \zeta_{\text{Pmin}} &  = \frac{1}{2} \left[1+\text{tanh}\left (\frac{S_{@\text{Pmin}}}{\rho_b} \right) \right ].
    \end{align}
    \item \textbf{Step 5}: the smoothed thrust and time rate of change of mass are written as 
    \begin{align}
        T^s_{\text{Ne@Pmax}} & =   N_{\text{e@Pmax}}^s \zeta_{\text{Pmax}} \frac{2 P_{\text{max}}  \eta_{\text{Pmax}}}{c_{\text{Pmax}}}, \\
        T^s_{\text{Ne@Pmin}} & =  N_{\text{e@Pmin}}^s \zeta_{\text{Pmin}} \frac{2  P_{\text{min}} \eta_{\text{Pmin}}}{c_{\text{Pmin}}},\\
 \dot{m}^s_{\text{Ne@Pmax}} & =  -\frac{T^s_{\text{Ne@Pmax}}}{c_{\text{Pmax}}}, \\
 \dot{m}^s_{\text{Ne@Pmin}} & =  -\frac{T^s_{\text{Ne@Pmin}}}{c_{\text{Pmin}}},
    \end{align}
    where $\eta_{\text{Pmax}}$ ($\eta_{\text{Pmin}}$) and $c_{\text{Pmax}}$ ($c_{\text{Pmin}}$) denote the efficiency and exhaust velocity values at maximum (minimum) power setting, respectively. Altogether, the thrust associated with maximum power is multiplied by two controllably smooth, activation coefficients in order to accomplish a smooth approximation of the number of engines (through $N_{\text{e@Pmax}}^s$) and a smooth approximation of the engine switching, $\zeta_{\text{Pmax}}$. The same procedure is followed for the thrust associated with minimum power. Once the thrust value is known, the time rate of change of mass is straightforward and is achieved by a simple division. Ultimately, these smooth approximation will be used to evaluate the total thrust and mass time rate of change as
    \begin{align}
    \textbf{T} & =  \hat{\bm{\alpha}} \left [ T^s_{\text{Ne@Pmax}} + T^s_{\text{Ne@Pmin}} \right ],\\
    \dot{m}  & = \dot{m}^s_{\text{Ne@Pmax}} + \dot{m}^s_{\text{Ne@Pmin}}.
\end{align}
\end{itemize}

At this stage, it is possible to determine all the information required for evaluating the RHS of Eq.~\eqref{eq:P2-statedynamics} and to propagate the differential equations to solve a TPBVP with smooth control inputs. The optimal direction of thrust steering vector is governed by the primer vector, which is a universal law when the direction of thrust vector is free and is determined as the definition of Lawden's primer vector:
\begin{align} \label{eq:pmvector}
    \hat{\bm{\alpha}} = - \frac{\bm{\lambda}^{\top} \mathbb{B}}{\|\bm{\lambda}^{\top} \mathbb{B} \|}.
\end{align}

In summary, the above steps are outlined to achieve smoothing at two levels: 1) a smooth transition between a finite set of operation modes, and 2) a smooth approximation of the switching function that commands engines to switch on or off. The combination of these two levels of smoothing determines the close neighbor of the optimal solution for the considered engine power settings. 

\subsection{Numerical Continuation Procedure}
Similar to the discussion in Section 4.4 of part 1 \cite{taheri2019anovelpart1}, there is a connection between the resulting smooth TPBVP and the original non-smooth OCP. Specifically, all of the control inputs are continuous and differentiable, which leads to having smooth, continuous state and costate dynamics. While these dynamics are rigorously smooth, the switches can be tuned by sweeping the smoothing parameter to be sufficiently sharp that they are qualitatively indistinguishable from the discontinuous controls that they approximate. In particular, we have a two-parameter family of smooth (superscript `s') neighboring OCPs that approach the actual OCP as the continuation parameters tend to zero. One can construct smoothed approximate neighbors that satisfy the invariant embedding constraint
\begin{align} \label{eq:connectionbetweensmooth}
    \begin{cases}
    \dot{\textbf{x}}^s = \dot{\textbf{x}}^s(\textbf{z}^s,\mathbf{U}^s,t)\\
    \dot{m}^s = \dot{m}^s(\textbf{z}^s,\mathbf{U}^s,t),\\
    \dot{\bm{\lambda}}^s = \dot{\bm{\lambda}}^s(\textbf{z}^s,\mathbf{U}^s,t),\\
    \dot{\lambda}_m^s = \dot{\lambda}_m^s(\textbf{z}^s,\mathbf{U}^s,t),
    \end{cases} \xrightarrow[\bm{\rho} = \bm{\rho}_{\text{max}} \rightarrow \mathbf{0}]{\mathbf{U}^s \rightarrow \mathbf{U}^* }
    \begin{cases}
    \dot{\textbf{x}}^* = \dot{\textbf{x}}^*(\textbf{z}^*,\mathbf{U}^*,t),\\
    \dot{m}^* = \dot{m}^*(\textbf{z}^*,\mathbf{U}^*,t),\\
    \dot{\bm{\lambda}}^* = \dot{\bm{\lambda}}^*(\textbf{z}^*,\mathbf{U}^*,t),\\
    \dot{\lambda}_m^* = \dot{\lambda}_m^*(\textbf{z}^*,\mathbf{U}^*,t),
    \end{cases}
\end{align}
where $\bm{\rho} = [\rho_b,\rho_c]^T$ denote the vector of smoothing parameters. Smooth and optimal control vectors (that consist of all control variables) are given as 
\begin{align*}
\mathbf{U}^s \in & \left [ \bm{\alpha}^s(\bm{\rho}), N_{\text{e@Pmax}}^s(\rho_b,\rho_c), N_{\text{e@Pmin}}^s(\rho_b,\rho_c)  \right ],\\
\mathbf{U}^* \in  & \left [ \bm{\alpha}^*(\bm{\rho}), N_{\text{e@Pmax}}^*(\rho_b,\rho_c), N_{\text{e@Pmin}}^*(\rho_b,\rho_c) \right ] \Bigg |_{\bm{\rho} = \bm{\rho}_{\text{min}} \approx \textbf{0}}.
\end{align*}

Note that the control vector consists of three components: direction of the thrust vector, $\bm{\alpha}$, number of engines operating at maximum power setting, $N_{\text{e@Pmax}}$, and number of engines operating at minimum power $N_{\text{e@Pmin}}$. As the continuation parameters approach zero, the smooth control vector, $\mathbf{U}^s$ approaches the values corresponding to the non-smooth OCP, $\mathbf{U}^*$. The numerical continuation procedure is straightforward and is analogous to methods widely used for solving optimal control problems \cite{mall2017epsilon,saghamanesh2018arobust}. In the next section, the results of applying the outlined set of steps on a low-thrust fuel-optimal trajectory optimization problem is demonstrated.

\subsection{Numerical Results for Earth-to-Comet 67P Problem: Same-Type Engine Cluster} \label{sec:resultssametype}
To demonstrate the utility of the proposed framework, a low-thrust trajectory from the Earth to comet 67P/Churyumov-Gerasimenko is considered. Due to the large change between the inclination, eccentricity, and semi-major axis orbital elements of the Earth and those of the comet, low-thrust trajectories would consist of more than one revolution around the Sun and may take up to 4 years \cite{ellison2018application}. Since we are going to take into account the perturbing acceleration due to all planets, the initial position vector cannot match that of the Earth. Therefore, it is assumed that the spacecraft is on the boundary of the SOI of the Earth with a positive along-the-track offset equal to one Earth SOI radius. Therefore, the initial position and velocity vectors of the spacecraft at $t = t_0$ are
\begin{align}
\textbf{r}_0 & = \begin{bmatrix}
-1671985.95664453\\
-151914424.309981\\
1699.37510504324
\end{bmatrix} \text{km}, & \textbf{v}_0 = \begin{bmatrix}
29.3070443053298\\
-0.596900982440449\\
-4.10911520334288 \times 10^{-4}
\end{bmatrix} \text{km/s}.
\end{align}

The final position and velocity vectors are
\begin{align}
\textbf{r}_f & = \begin{bmatrix}
-465627493.14461\\
-50530561.3073027	\\
40190127.9500019
\end{bmatrix} \text{km}, & \textbf{v}_f = \begin{bmatrix}
-9.7217789589445\\
-14.6294809300934\\
-0.234945260833124
\end{bmatrix} \text{km/s}.
\end{align}

The spacecraft is assumed to leave the Earth on June 17, 2024, and would reach the comet on April 28, 2029. The time of flight is fixed at $t_f - t_0 = 1776$ days. The initial mass of the spacecraft is $m_0 = 3000$ kg, the beginning-of-life power is set to $P_{0,\text{BOL}} = 30$ kW and $\sigma = 2$\% per year. It is assumed that 500 Watts of power is used to energize the PPU and operate various spacecraft sub-systems during the entire time of flight. For the BPT-4000 Extra-High-$I_{\text{sp}}$ engine, the following values are used: $\eta_{\text{Pmax}} = 0.48576$, $\eta_{\text{Pmin}} = 0.14806$, $c_{\text{Pmax}} = 25491.663$ (m/s), $c_{\text{Pmin}} = 6158.059$ (m/s), $P_{\text{max}} = 4839$ Watts, $P_{\text{min}} = 302$ Watts. These parameters have been used to generate the power levels summarized in Table \ref{tab:E67p_Cases}.

Planetary perturbations in the modeling represents the disturbing acceleration due to all of the planets of the Solar System from the innermost planet Mercury to the outermost planet Neptune, $N_{\text{sb}} = 8$. The numerical CX method (with its implementation outlined in part 1) is used to evaluate the contribution of planetary perturbations into costate dynamics.

\begin{table}[htbp!] 
\begin{center} 
		\caption{Summary of the results for the Earth-to-67P problem with a cluster of four BPT-4000 Extra-High-$I_{\text{sp}}$ thrusters; $\rho_{\text{b}} = \rho_{\text{c}}= 1.0 \times 10^{-2}$. } \label{tab:E67p_Cases}
		{\small
		\begin{tabular}{c c c c c c}
        \hline
        \hline
         Case  & Two-body & Power  & Degradation  & Planetary    & $m_f$\\
          \#       &          & Model  & Model        & Perturbation & (kg) \\    
                 &  $\mu_{\text{sun}}\textbf{r}/r^3 $        &  $\phi(r) = 1/r^2$  &  $\psi(t)$        & $\textbf{a}_{\text{sb}}$ & \\
         \hline          
          1      & Yes      & Yes         & No                 & No  & 1930.507  \\
          2      & Yes      & Yes         & Yes                & No  & 1922.064 \\
          3      & Yes      & Yes         & Yes                & Yes & 1922.301 \\
        \hline
        \hline
        \end{tabular}
		}
	\end{center}
\end{table}

To quantify the impact of various models, and given the flexibility of the tool, three cases are considered and are listed in Table \ref{tab:E67p_Cases}. The difference in these cases is due to the inclusion of solar power system degradation and planetary perturbation models summarized as follows
\begin{itemize}
    \item \textbf{Case 1}: two-body gravitational model without  consideration of variation of power due to change in distance and degradation of the solar arrays, and no inclusion of planetary perturbations,
    \item \textbf{Case 2}: two-body gravitational model with consideration of variation of power due to change in distance and degradation of the solar arrays, and no inclusion of planetary perturbations,
    \item \textbf{Case 3}: two-body gravitational model with consideration of variation of power due to change in distance and degradation of the solar arrays, and with inclusion of planetary perturbations.
\end{itemize}

\begin{figure}[htbp!]
\centering
\includegraphics[width=4.0in]{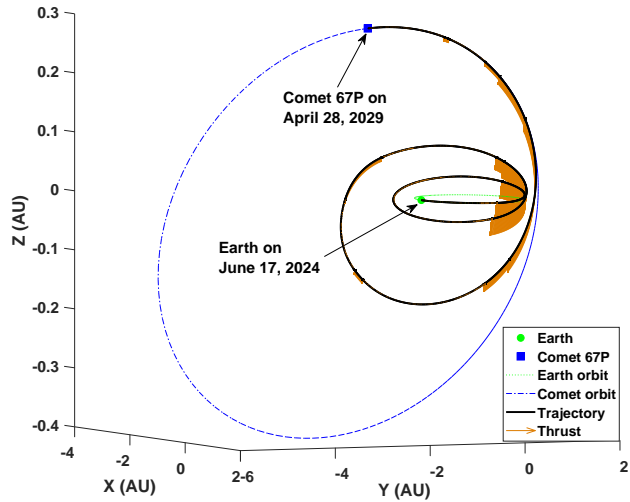}
\caption{Same-type engines: trajectory from the Earth to comet 67P case 3 with 4 $\times$  BPT-4000 Extra High-$I_{\text{sp}}$ engines.}
\label{fig:E67PNe4Case3Traj}
\end{figure}

Figure \ref{fig:E67PNe4Case3Traj} depicts the location of the Earth on its orbit at the time of departure (June 17, 2024), low-thrust trajectory, and location of the comet on its orbit at the end of flight (April 28, 2029), all in the heliocentric J2000 frame of reference. The optimal solution corresponds to making two revolutions around the Sun, i.e., $N_{\text{rev}} = 2$. 

\begin{figure}[htbp!]
\centering
\includegraphics[width=3.8in]{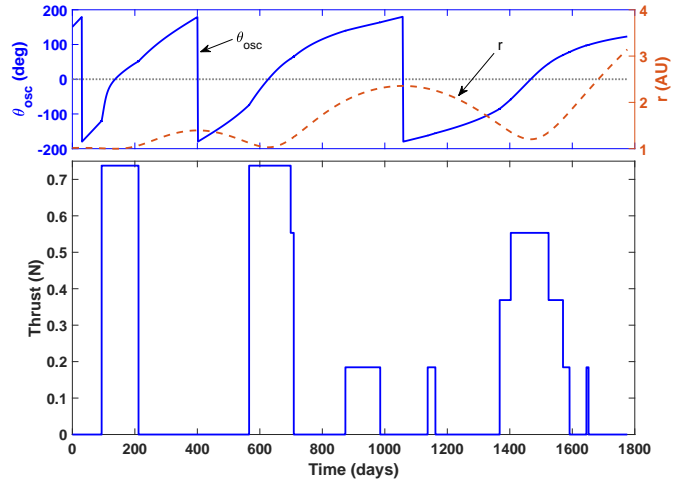}
\caption{Same-type engines: time history of osculating true anomaly and thrust for case 3 with 4 $\times$  BPT-4000 Extra High-$I_{\text{sp}}$ engines.}
\label{fig:E67PNe4Case3TA_Thrust}
\end{figure} 

\begin{figure}[htbp!]
\centering
\includegraphics[width=4.0in]{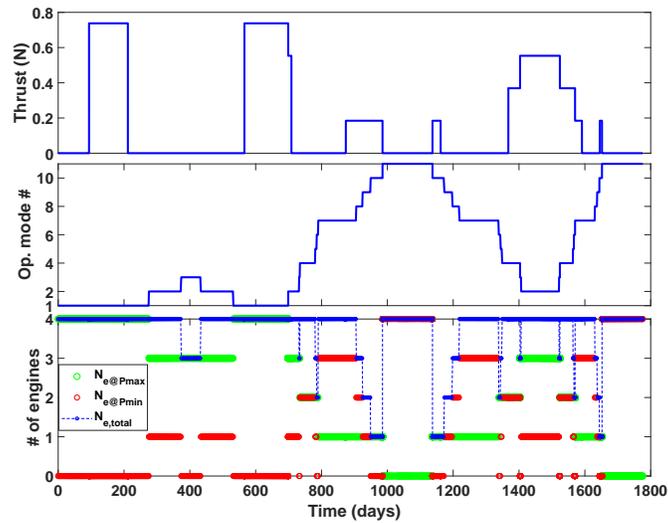}
\caption{Same-type engines, time history of the power levels and number of potential engines on the optimal trajectory for case 3 with 4 $\times$  BPT-4000 Extra High-$I_{\text{sp}}$ engines.}
\label{fig:E67PNe4Case3Powerlevel_Nengine}
\end{figure} 
Figure \ref{fig:E67PNe4Case3TA_Thrust} shows the time history of the osculating true anomaly and the thrust profile for the optimal trajectory. The majority of thrusting occurs at the perihelion passages. Figure \ref{fig:E67PNe4Case3Powerlevel_Nengine} shows the time history of the operation modes and the number of engines. Only the first 11 operation modes are engaged during this ``optimal'' trajectory. We draw your attention to the fact that the trajectory consists of zero-thrust arcs at the beginning and at the end of the trajectory. This indicates that the particular configuration of the propulsion system is more capable and the trajectory starts and terminates on the so-called \textit{late-departure} and \textit{early-arrival} boundaries. These boundaries will trace a curve and are revealed as part of the optimal switching surfaces introduced by Taheri and Junkins \cite{taheri2018how}. 

\begin{figure}[htbp!]
\centering
\includegraphics[width=3.5in]{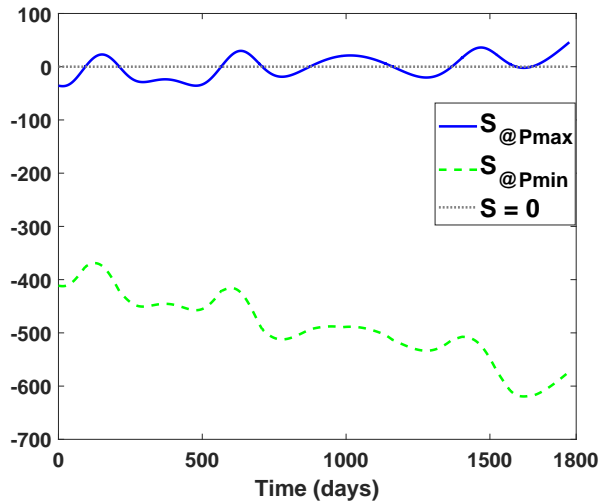}
\caption{Same-type engines: time history of the switching functions for case 3 with 4 $\times$  BPT-4000 Extra High-$I_{\text{sp}}$ engines.}
\label{fig:E67PNe4Case3SwitchingFunctions}
\end{figure}

Another important point is that the particular number of engines associated with each operation mode is just an indication of the number and operation power setting of \textit{potential} engines. It does not necessarily mean that all of the engines for that particular operation mode will become active. The activation depends on optimality criteria and depends on the sign of the respective switching functions. In order to clarify this point, Figure \ref{fig:E67PNe4Case3SwitchingFunctions} depicts the time history of the switching functions defined in Eqs.~\eqref{eq:sfpmax} and \eqref{eq:sfpmin}. In fact, it shows that only those engines that operate at the maximum power are contributing to the thrust and rate of change of mass, whereas the engine(s) at the lowest power setting is (are) always inactive. Solar arrays may have generated enough power to turn engines on, however, the optimality principle sets the command to turn an engine (or a set of engines) on or off. 

At the end of trajectory, $S_{\text{@Pmax}}$ is still positive. However, the $11^{\text{th}}$ operation mode is engaged (see middle sub-plot in Figure \ref{fig:E67PNe4Case3Powerlevel_Nengine}), which corresponds to having only 4 engines that can only operate at $P_{\text{min}}$ (see Table \ref{tab:representativepowerlevels}). The combination of the mentioned factors leads to a zero-thrust arc.

\subsection{Same-Type Engine Clusters With Large Number of Engines}
In order to test the performance of the CSC framework, we studied cases in which the number of engines was increased to 6, 10 and a hypothetical cluster with 20 engines. A condensed list of operation modes of each cluster is given in Table \ref{tab:powermodeswithlargerNe}. The total number of operation modes can be further reduced if the feasible range of distance from the Sun is taken into account. Here we have considered the total number of operation modes purely based on the combination of engines.
\begin{table}[htbp!] 
\begin{center} 
		\caption{Representative first and last operation modes for a trajectory from the Earth to comet 67P with multiple BPT-4000 Extra-High-$I_{\text{sp}}$ engines ($N_e \in \{6,10, 20\}$) and the optimal final mass  for case 1; power unit is Watts. } \label{tab:powermodeswithlargerNe}
		{\small
		\begin{tabular}{l c c c c | c}
        \hline
        \hline
         $N_e$ & Mode \#         & $P_{\text{used}} $  &  $N_{\text{e@Pmax}}$  &  $N_{\text{e@Pmin}}$  & $m_f$ (kg)\\
         \hline          
        \multirow{2}{*}{6} & 1 & 29,034 & 	6 &	0 & \multirow{2}{*}{1939.844} \\
                          & 27 & 302 &	0 &	1 & \\
         \hline
         \multirow{2}{*}{10}& 1 & 43,853	 & 9 &	1 & \multirow{2}{*}{1939.844}\\
                           & 64 & 302 &	0 &	1   &  \\
         \hline
        \multirow{2}{*}{20} & 1 & 46,873 &	9 &	11 &   \multirow{2}{*}{1939.815} \\
                           & 164 & 302 &	0 &	1  & \\
         
        \hline
        \hline
        \end{tabular}
		}
	\end{center}
\end{table}

Of course, having 20 engines is likely not realistic, but we use this scenario as a verification and demonstration step to assess the performance of the CSC framework for such a problem with 164 different operation modes. No other known optimal control approach is known that can accommodate these many propulsion system operation modes. The dry mass of each engine is about 7.5 kg, which adds up to 150 kg if a vehicle design analysis is conducted. More importantly, the required power would impact the dimension and weight of the solar arrays. However, the same vehicle/engine parameters of the previous section are used. In order to simplify the results, only case 1, is considered where the solar array degradation and third-body perturbations are ignored. The last column of Table \ref{tab:powermodeswithlargerNe} presents the final mass.

\begin{figure}[htbp!]
\centering
\includegraphics[width=4.0in]{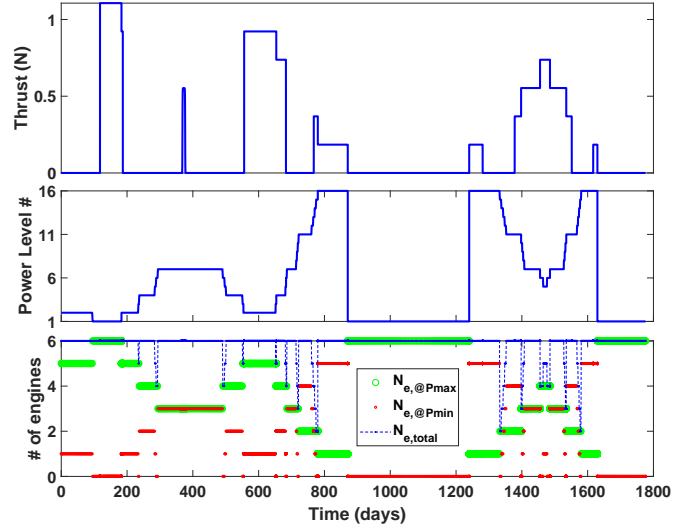}
\caption{Same-type engines: the thrust, power levels and number of potential operating engines vs. time for case 1 with 6 $\times$ BPT-4000 Extra High-$I_{\text{sp}}$ engines.}
\label{fig:E67PNe6Case1ThrustPowerlevel_Nengine}
\end{figure}

\begin{figure}[htbp!]
\centering
\includegraphics[width=4.0in]{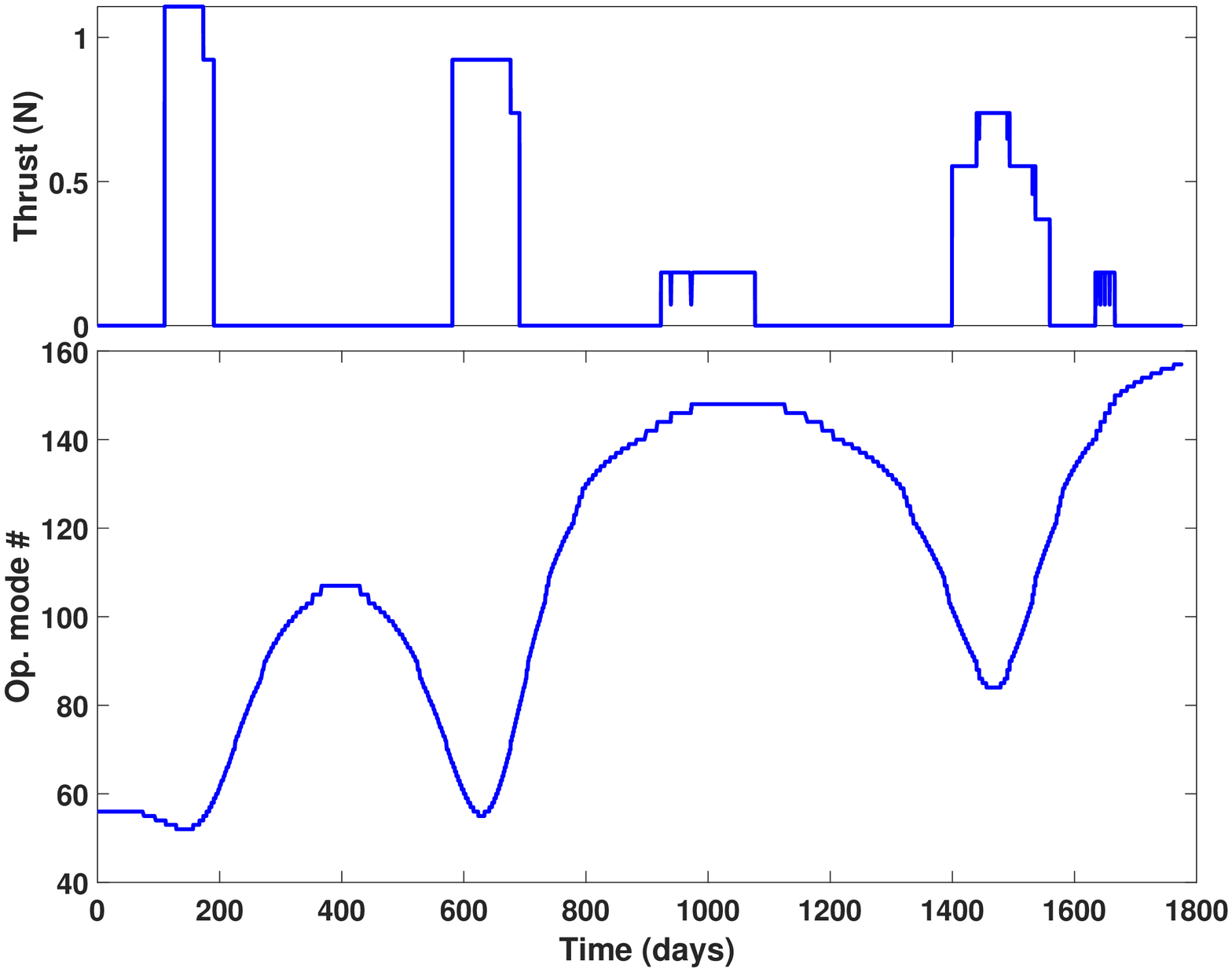}
\caption{Same-type engines: the thrust, power levels, and number of potential operating engines vs. time for case 1 with 20 $\times$ BPT-4000 Extra High-$I_{\text{sp}}$ engines engines.}
\label{fig:E67PNe20Case1ThrustPowerlevel_Nengine}
\end{figure}

Figures \ref{fig:E67PNe6Case1ThrustPowerlevel_Nengine} and \ref{fig:E67PNe20Case1ThrustPowerlevel_Nengine} depict the time history of thrust, power levels, and number of engines for solutions with 6 and 20 engines. Figure \ref{fig:E67PNe6Case1ThrustPowerlevel_Nengine} shows that when $N_e = 6$, only the first 22 (out of the total 27) operation modes are engaged. However, when $N_e = 20$, operation modes \#52 through \#157 are engaged. The thrust profiles in Figures \ref{fig:E67PNe6Case1ThrustPowerlevel_Nengine} and \ref{fig:E67PNe20Case1ThrustPowerlevel_Nengine} show similar trends especially during the first thrust arcs and up to the last two the thrust arcs. While a complicated set of switches is evident (the third sub-figure of Figure \ref{fig:E67PNe20Case1ThrustPowerlevel_Nengine}), the resulting optimal 20-engine thrust profile is qualitatively very similar to the 6-engine profile of Figure \ref{fig:E67PNe6Case1ThrustPowerlevel_Nengine}.

Of course, in our analysis, the inert mass variations due to having a greater number of engines are not taken into account. Therefore, the final masses reported in Table \ref{tab:powermodeswithlargerNe} do not differ from each other significantly (and are comparable to Case 1 in Table \ref{tab:E67p_Cases}). This demonstration is simply used to verify the applicability of our approach. The results indicate that the proposed smoothing method is capable of handling discrete set of power modes as high as 164. The time history of thrust for a cluster with 10 engines (with 64 operating modes) is identical to the one in Figure \ref{fig:E67PNe6Case1ThrustPowerlevel_Nengine} and is not plotted. The switching functions also are identical to those depicted in Figure \ref{fig:E67PNe4Case3SwitchingFunctions}.


\section{Mixed/Different-Type Engine Clusters} \label{sec:differentype}
A more complicated strategy is needed if/when all of the engines within a cluster are not of the same type. In these situations, each individual engine has to be taken into account twice: one with its power input set at $P_{\text{min}}$ and one with its power input set to $P_{\text{max}}$. Each of these two settings are multiplied by a smoothed activation function, $\gamma$ (with appropriate indexing). As a consequence, Eq.~\eqref{eq:sumthrust} can be written as
\begin{align}
    \textbf{T} & = \hat{\bm{\alpha}}  \sum_{i=1}^{Ne} \left [ \gamma_{\text{i@Pmin}} T_{\text{i@Pmin}} + \gamma_{\text{i@Pmax}} T_{\text{i@Pmax}} \right ], \label{eq:thrustmixed}\\    
    \dot{m} & = ~~\sum_{i=1}^{Ne} \left [ \gamma_{\text{i@Pmin}} \dot{m}_{\text{i@Pmin}} + \gamma_{\text{i@Pmax}} \dot{m}_{\text{i@Pmax}} \right ] \label{eq:mdotmixed} 
\end{align}
where the unknown control inputs are $\hat{\bm{\alpha}}$, and $\gamma_{\text{i@Pmin}}$ and $\gamma_{\text{i@Pmax}}$ (for $i = 1, \cdots, N_e$).

\subsection{Implementation Details for Mixed-Type Engine Clusters} \label{sec:mixethrustermethodology}
Unlike a cluster with the same type of engines, we need to determine the activation function of individual engines. In addition, each engine has two distinct power settings (excluding $P_{\text{en}} = 0$). Therefore, in general we can not group them into two categories and turn the problem into the task of determining the number of active engines in that particular category. Here we need to deal with each engine's on/off switching individually. 

Immediately after the activation functions are known, we can substitute their values into Eqs.~\eqref{eq:thrustmixed} and \eqref{eq:mdotmixed} and evaluate the RHS of the set of state-costate dynamics and propagate them numerically to solve the shooting problem. Again, we choose a strategy to cast the problem into a set of operation modes. The only difference is that the columns of the operating table will grow depending on the number of considered engines.    

Since we are studying engine clusters with a mixed sub-cluster or distinct number of engines, we specify each cluster with an ordered list of engine identifiers, IDs, that are defined in Table \ref{tab:enginepar}. For instance, an engine sequence $E = \{1,2,3,4\}$ denotes an engine cluster that consists of the three variants of BPT-4000 engines, i.e., High-$I_{\text{sp}}$, High-Thrust, Extra High-$I_{\text{sp}}$, and one Next TT10 High-$I_{\text{sp}}$ engine. We can also have an engine cluster $E = \{3,4,6\}$, which consists of three distinct types of engines.

The logic for generation of table of operation modes is not as straightforward as it is for the same-type engine clusters. In particular, it turns out that for certain engine sequences, it is possible to have different combinations of engines operating at different power settings, but require identical $P_{\text{used}}$. They consist of different power settings for engines, but the aggregate used power is the same. Such a situation can occur, for example, if we have a cluster of two BPT-4000 engines with a sequence $E = \{2,3\}$. Despite the fact that these two engines have different characteristic performance curves, they operate over almost the same power range, i.e., $P_{\text{en}} \in [0.302, 4.839]$, which means that they have identical maximum and minimum power settings. 

Table \ref{tab:powerlevelsconflict} summarizes the original data table for $E = \{2,3\}$ where three different sets of operation modes are distinguishable with the same value of $P_{\text{used}}$. There are four columns since there are two engines (each can operate at maximum or minimum power). The first operation mode, for example, corresponds to when both engines operate at their maximum power settings, $P_{\text{max}}$. Note that the presented table is just a representative example and there could be situations in which multiple operation modes have the same value of $P_{\text{used}}$ if there are more engines. In the end, we need to have a strictly monotonic (decreasing/increasing with respect to power) list of operation modes, where the power of each mode ($P_{\text{used}}$) is different from the other operation modes. Note two distinct modes consuming the same power do not consume the same fuel mass. Since we are interested in fuel-optimal trajectories, we break the tie between the operation modes (among the conflicting sets) by selecting the operation mode with the least value of the rate of change of mass, $\dot{m}$. In the given example, operation modes \# 3, 4, and 8 will pass the fuel optimality criterion filtering process. Table \ref{tab:powerlevelscpostfilter} summarizes the resulting operation modes for $E = \{2,3\}$.

\begin{table}[htbp!] 
\begin{center} 
		\caption{Representative table of operation modes for a trajectory from the Earth to comet 67P with engine sequence $E = \{2,3\}$ and when some of the modes are conflicting. } \label{tab:powerlevelsconflict}
		{\small
		\begin{tabular}{c c c c c c c }
        \hline
        \hline
         \multirow{3}{*}{Mode \#}   & \multirow{3}{*}{$P_{\text{used}}$ (Watts)}  &  \multicolumn{4}{c}{Engine ID}     &  $\dot{m}$ \\
         \cline{3-6} 
              &                     &  \multicolumn{2}{c}{2}    & \multicolumn{2}{c}{3}  &(mg/s) \\
         \cline{3-4} \cline{5-6}
                &                     & $P_{\text{min}}$ & $P_{\text{max}}$ & $P_{\text{min}}$ & $P_{\text{max}}$ & \\
         \hline          
          1 & 9,678   & 0 & 1 &	0 &	1 & 22.59 \\
           \rowcolor{LightCyan}2 & 5,141   & 0 &	1 & 1 &	0 & 17.72  \\
           \rowcolor{LightCyan}3 & 5,141   & 1 &	0 &	0 &	1 & 9.51  \\
          \rowcolor{Gray}4 & 4,839   & 0 &	0 &	0 &	1 & 7.23  \\
          \rowcolor{Gray}5 & 4,839   & 0 &	1 &	0 &	0 & 15.36  \\
          6 & 604     & 1 & 0 & 1 &	0 & 4.64  \\
          \rowcolor{LightCyan}7 & 302     & 0 &	0 &	1 &	0 & 2.36  \\
          \rowcolor{LightCyan}8 & 302     & 1 &	0 &	0 &	0 & 2.28  \\
        \hline
        \hline
        \end{tabular}
		}
	\end{center}
\end{table}

\begin{table}[htbp!] 
\begin{center} 
		\caption{Post-filtering table of operation modes for a trajectory from the Earth to comet 67P with engine sequence $E = \{2,3\}$. } \label{tab:powerlevelscpostfilter}
		{\small
		\begin{tabular}{c c c c c c c }
        \hline
        \hline
         \multirow{3}{*}{Mode \#}   & \multirow{3}{*}{$P_{\text{used}}$ (Watts)}  &  \multicolumn{4}{c}{Engine ID}     &  $\dot{m}$ \\
         \cline{3-6} 
              &                     &  \multicolumn{2}{c}{2}    & \multicolumn{2}{c}{3}  &(mg/s) \\
         \cline{3-4} \cline{5-6}
                &                     & $P_{\text{min}}$ & $P_{\text{max}}$ & $P_{\text{min}}$ & $P_{\text{max}}$ & \\
         \hline          
          1 & 9,678   & 0 & 1 &	0 &	1 & 22.59 \\
          2 & 5,141   & 1 &	0 &	0 &	1 & 9.51  \\
          3 & 4,839   & 0 &	0 &	0 &	1 & 7.23  \\
          4 & 604     & 1 & 0 & 1 &	0 & 4.64  \\
          5 & 302     & 1 &	0 &	0 &	0 & 2.28  \\
        \hline
        \hline
        \end{tabular}
		}
	\end{center}
\end{table}

As far as the implementation is concerned, the same overall steps outlined in Section \ref{sec:impdetails} are followed to compute the thrust and time rate of change of mass. The only major difference is that the operations are extended over a larger set of data points (as opposed to having only two category of engines). Specifically, element-wise multiplications of matrices and summation over columns of the resulting operations are performed. Note that $\textbf{w}$ is a matrix obtained by repeating its vector form $k$ times (where $k = 2 N_e$) to construct a $\sigma \times k$ matrix. The time rate of change of mass and thrust associated with individual engine power settings have to be calculated. The details are removed for brevity and it is straightforward to establish the required relations. The ultimate relation can be written as
\begin{align}
    \textbf{T} & = \hat{\bm{\alpha}}  \sum_{i=1}^{Ne} \left [ \gamma_{\text{i@Pmin}} \zeta_{\text{i@Pmin}} T_{\text{i@Pmin}} + \gamma_{\text{i@Pmax}} \zeta_{\text{i@Pmax}} T_{\text{i@Pmax}} \right ], \\    
    \dot{m} & = ~~\sum_{i=1}^{Ne} \left [ \gamma_{\text{i@Pmin}} \zeta_{\text{i@Pmin}} \dot{m}_{\text{i@Pmin}} + \gamma_{\text{i@Pmax}}  \zeta_{\text{i@Pmax}} \dot{m}_{\text{i@Pmax}} \right ], 
\end{align}
where each engine power setting has its own coefficient triggered also by the switching function for that particular engine and power setting, which are the $\zeta$ coefficients in the above relations. The direction of the total thrust vector is still governed by Eq.~\eqref{eq:pmvector}. As far as the invariant embedding is concerned, Eq.~\eqref{eq:connectionbetweensmooth} still holds and we use a two-parameter family of OCPs. The only difference is that the control vector consists of additional elements given as
\begin{align*}
\mathbf{U}^s = 
\begin{bmatrix} 
 \bm{\alpha}^s \\
 \gamma_{\text{1@Pmin}}^s\\ 
 \gamma_{\text{1@Pmax}}^s\\
 \vdots \\
 \gamma_{\text{N}_e\text{@Pmin}}^s\\ 
 \gamma_{\text{N}_e\text{@Pmax}}^s\\
\end{bmatrix}_{\bm{\rho}} \xrightarrow{\bm{\rho} \rightarrow \bm{\rho}_{\text{min}}} \mathbf{U}^* = 
\begin{bmatrix} 
 \bm{\alpha}^* \\
 \gamma_{\text{1@Pmin}}^*\\ 
 \gamma_{\text{1@Pmax}}^*\\
 \vdots \\
 \gamma_{\text{N}_e\text{@Pmin}}^*\\ 
 \gamma_{\text{N}_e\text{@Pmax}}^*\\
\end{bmatrix}_{\bm{\rho} = \bm{\rho}_{\text{min}} \approx \textbf{0}},
\end{align*}
where the $\gamma$ values in the control vectors have absorbed the $\zeta$ values.

\subsection{Numerical Results for Earth-to-Comet 67P Problem: Different-Type Engine Clusters}
The methodology outlined in Section \ref{sec:mixethrustermethodology} is used for designing trajectories from the Earth to comet 67-P. The boundary conditions of the problem, time of flight, and the initial mass of the spacecraft are all identical to those in Section \ref{sec:resultssametype}. 

Table \ref{tab:difftypesummary} summarizes the results for eight clusters with different engines (lowest 2 to largest 4 engines) and with different combinations of engines. A subscript index number is assigned to each engine sequence. For each sequence, the total number of operation modes, the number of engaged modes, and the final mass is reported. Recall also that engagement of a mode means that those engines are the best combination, but activation of each engine is triggered by its switching function. The result of the engine sequence, $E_6$, is identical to the same-engine result, but is obtained using the proposed scheme for a different-type engine cluster.  
\begin{table}[htbp!] 
\begin{center} 
		\caption{Summary of different-type engine clusters the Earth-to-comet 67P problem case 1. } \label{tab:difftypesummary}
		{\small
		\begin{tabular}{l c c c c }
        \hline
        \hline
         Cluster sequence   & Total \# of   &  Range of  &  $m_f$ \\
                    & Opt. Modes  &   Engaged Modes &  (kg) \\
         \hline
         $E_{1} =\{2,3\}$     & 5    & 1-4  & 1726.413 \\
         $E_{2} =\{3,5\}$     & 8    & 1-6  & 2029.389 \\
         $E_{3} =\{4,5\}$     & 8    & 1-6  & 2152.619 \\
         $E_{4} =\{2,4,5\}$        & 26   & 1-20 & 2173.261 \\
         $E_{5} =\{3,4,5\}$   & 26   & 1-20 & 2192.719 \\
         $E_{6} =\{3,3,3,3\}$ & 14   & 1-11 & 1930.507  \\
         $E_{7} =\{2,2,3,3\}$ & 14   & 1-11 & 1844.453	\\
         $E_{8} =\{2,3,4,5\}$ & 53  & 1-43 & 2159.789  \\
         \hline
        \hline
        \end{tabular}
		}
	\end{center}
\end{table}

\begin{figure}[htbp!]
\centering
\includegraphics[width=5.0in]{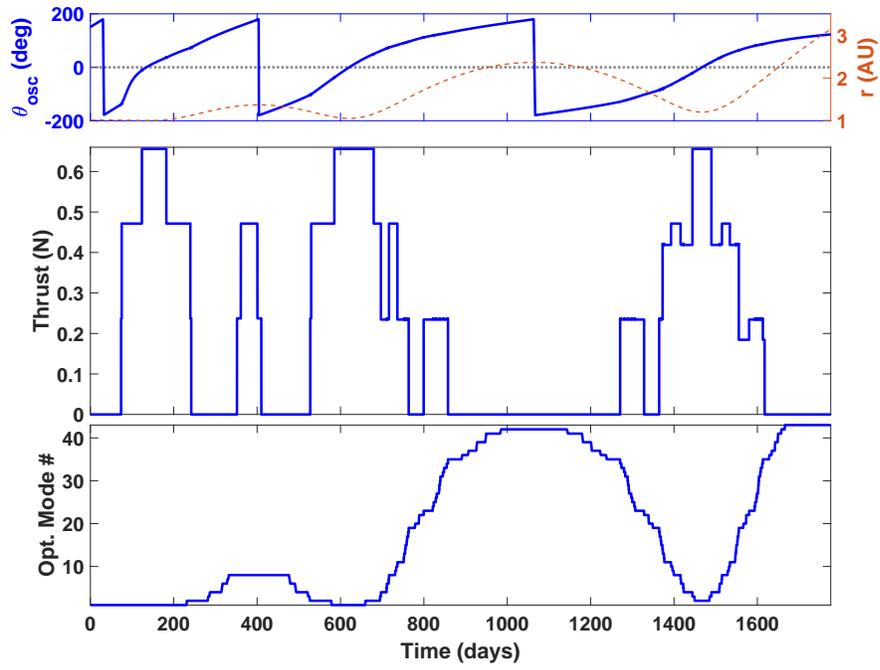}
\caption{Different-type engine cluster: time history of osculating true anomaly, thrust, and operation modes for engine sequence $E_8 = \{2,3,4,5\}$ case 1 with $\rho_b = \rho_c = 0.01$.}
\label{fig:E67PSeq8Case1TA_Thrust_Opmodesfig}
\end{figure}

\begin{figure}[htbp!]
\centering
\includegraphics[width=5.0in]{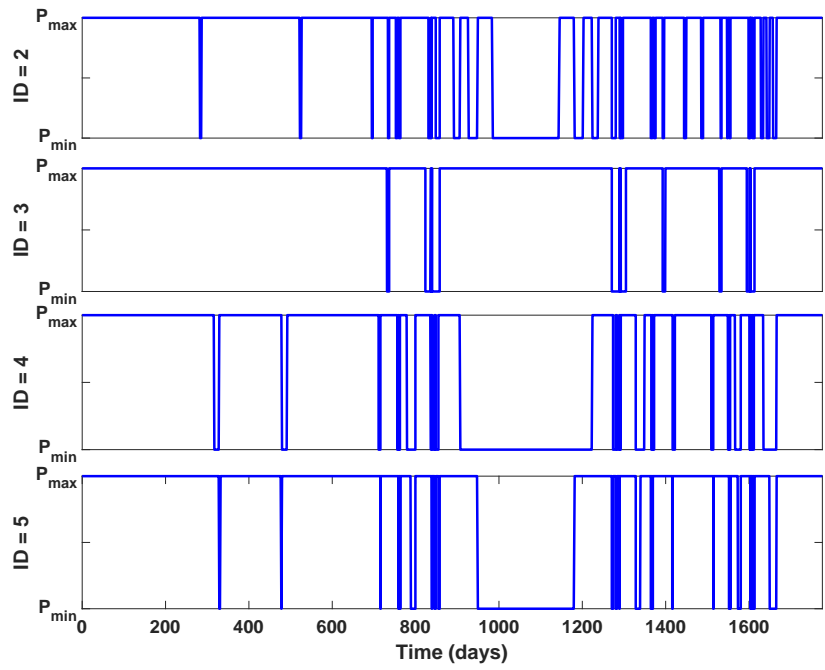}
\caption{Different-type engine cluster: switches between the engine power modes as the transitions occur between operation modes for cluster sequence $E_8 = \{2,3,4,5\}$ case 1.}
\label{fig:E67PSeq8Case1EngineSwitches}
\end{figure}

\begin{figure}[htbp!]
\centering
\includegraphics[width=4.0in]{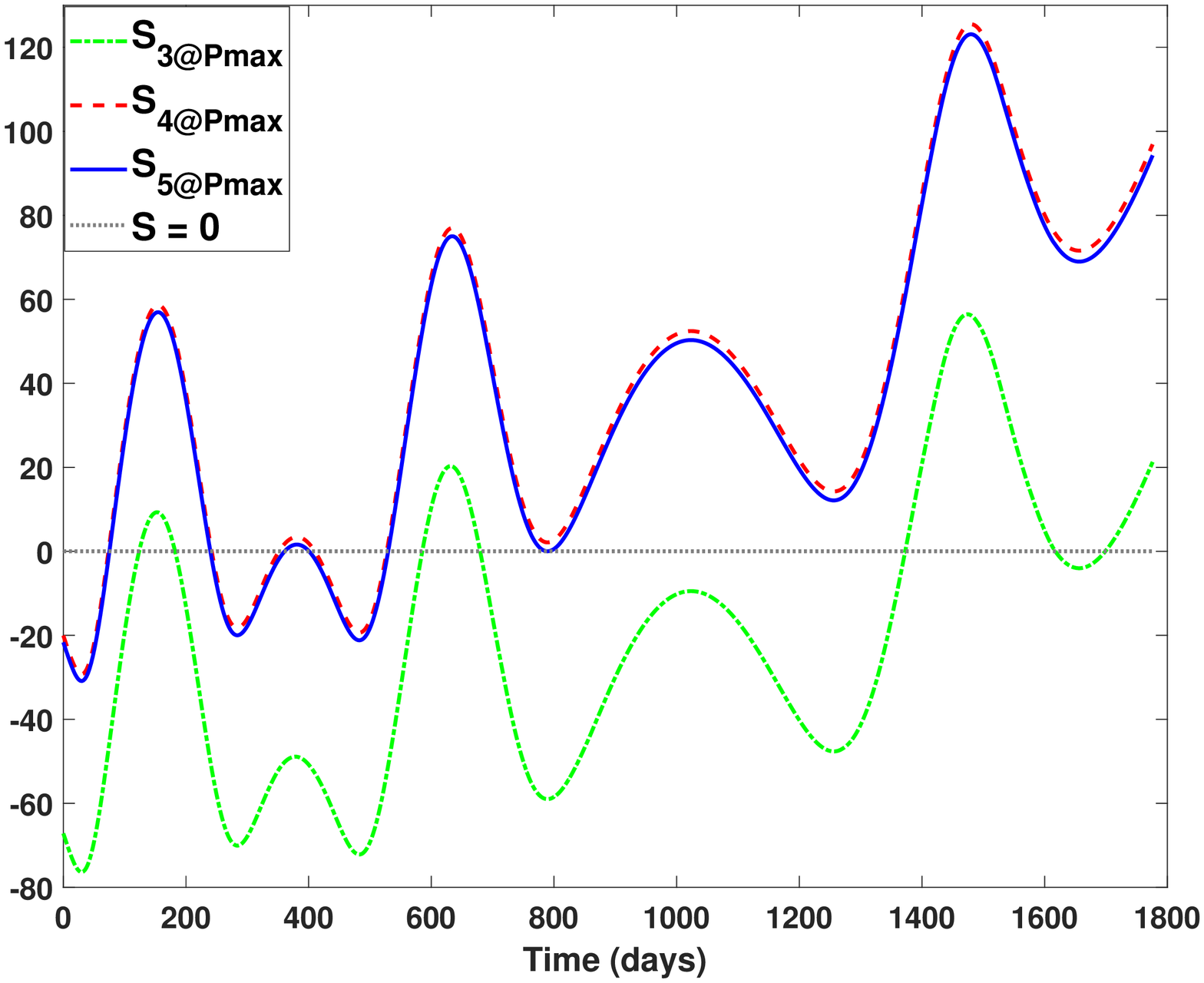}
\caption{Different-type engine cluster: time history of the switching functions that contribute to the solution for cluster sequence $E_8 = \{2,3,4,5\}$ case 1.}
\label{fig:E67PSeq8Case1SFs}
\end{figure}

\begin{figure}[htbp!]
\centering
\includegraphics[width=5.0in]{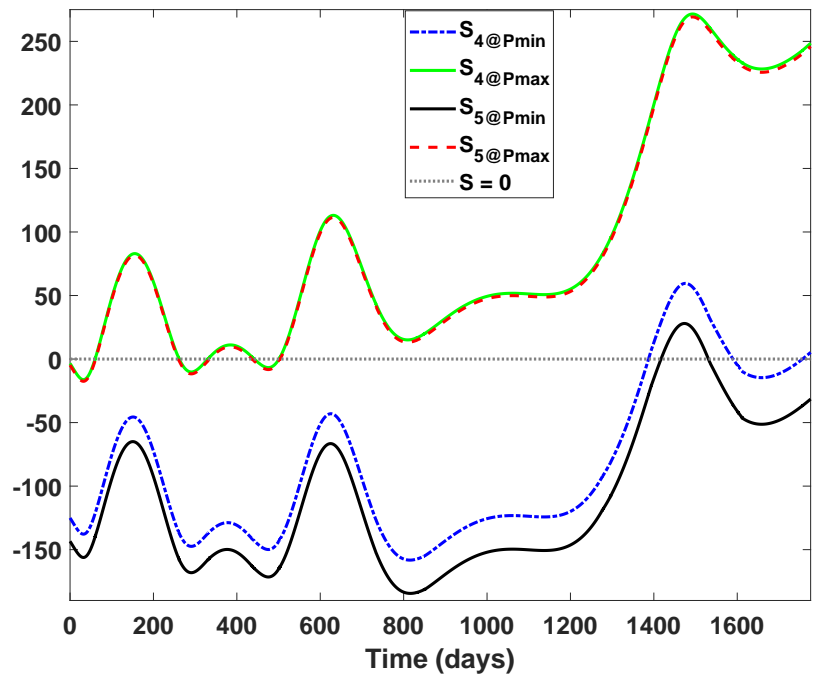}
\caption{Different-type engine cluster: time history of the switching functions for cluster sequence $E_3 = \{4,5\}$ case 1.}
\label{fig:E67PSeq3Case1SFs}
\end{figure}

Figure \ref{fig:E67PSeq8Case1TA_Thrust_Opmodesfig} shows the details of the trajectory for the engine sequence, $E_8$ and Figure \ref{fig:E67PSeq8Case1EngineSwitches} shows the switches between the maximum and minimum power settings corresponding to the total engaged operation modes. Figure \ref{fig:E67PSeq8Case1SFs} depicts the time history of the switching functions of those engines that contribute to the trajectory. This sequence has a total number of eight switching functions. Once again, only engines with IDs 3, 4 and 5 become active, which means that we could have removed engine ID 2 from this particular engine cluster. This type of analysis is quite instrumental for configuring a suitable engine cluster. 
In a majority of the cases engines operate at their maximum power settings. However, it is also possible for some of the engines to operate at their minimum power settings. For instance, when the engine sequence $E_3$ is used, the time history of the switching functions indicates that all of the engine operate at their maximum and minimum power settings. In particular, the minimum power settings of both engines are activated during the final time duration as is shown in Figure \ref{fig:E67PSeq3Case1SFs}. 

Obviously, we can not have an engine operating at its maximum and minimum power settings simultaneously. Recall that it is the combined effects of the switching function and the engaged operation mode that determines the actual operation of an engine. This is identical to the analysis presented for the zero-thrust arc in Figure \ref{fig:E67PNe4Case3SwitchingFunctions}, where the switching function was not the sole indicator of operation.

Engine IDs 4 and 5 correspond to Next TT10 and NEXT TT11 engines that have a large separation between their power settings compared to the other engines, which is one of the reasons that minimum power modes appear in solutions associated with this particular engine cluster. 

While restricting the admissible control to minimum or maximum power setting is not guaranteed to be the optimal strategy (for these throttleable engines), it may be viewed as a first effort to solve such challenging problems. In fact, we can add or remove the number of power settings of each engine at our will (from its current two options, $P_{\text{en}} \in \{P_{\text{min}}, P_{\text{max}}\}$ to beyond two) so long as the number of total operation modes do not grow to beyond some computational limit (since we are using a single-shooting scheme). It is also possible to set the value for the upper and lower power limits if we are interested to perform any particular type of analysis. For instance, if we set $P_{\text{min}} = 4.0$ KWatts instead of its original value, $P_{\text{min}} = 0.302$ KWatts, the final mass of the engine sequence $E = \{3,5\}$ becomes $m_f = 2031.43$ kg, which is 2.04 kg (0.1\%) greater than the value reported in Table \ref{tab:difftypesummary}. 

There are likely better strategies, but the proposed strategy and results of this study are still relevant and of importance, establish a beginning framework to impose a larger set of admissible control modes for multi-engine clusters. In particular, it serves as a first effort to solve such problems in the indirect optimization methods. The results can be used to gain insights about the structure of the control, which is helpful for formulating and solving the actual optimal control problem.

Our future work is focused on alternative/improved approaches for finding solutions that are closer to ``optimal'' for the case that multiple engines can be individually throttled. We should mention, however, that with respect to the selected discrete modes of operation of engines, the solutions are at least approximate local extremals since with sufficiently small $\rho_b$ and $\rho_c$, they satisfy the first-order necessary conditions of optimality.

%% file: Conclusion.tex
\section{Conclusion} \label{sec:conclusion}
In this work, the problem of designing spacecraft trajectories that exploit the indirect method of optimal control and is capable of handling multiple discrete operating modes is studied. In particular, trajectory optimization problems are considered for spacecraft with a cluster of engines for which the number of operating engines and the power at which each engine operates are not known \textit{a priori}. 

Application of the Composite Smooth Control (CSC) framework is presented that enables reducing the problem to a two-point boundary-value problem and facilitates the numerical solution of this problem. A two-parameter family of smoothly switching mode accommodates the case of multiple identical or discrete engine design, such that a very close neighbor of the Pontryagin necessary conditions, and the associated boundary conditions, can be efficiently determined. In the case of a spacecraft with a cluster of engines, a pre-processing step was proposed to define distinct modes of operation and then apply the CSC framework to determine the number and operating power of the engines. This framework is evidently the first practical method for applying indirect optimal control approach to multi-mode systems of this level of complexity, with realistic polynomial engine performance curves, and when a single-shooting scheme is used.  

The results indicate that the proposed framework performs well for solving problems in which the spacecraft is equipped with clusters consisting of 4, 6 or 10 engines. Application of the method is also demonstrated for an extreme hypothetical case of a cluster of 20 engines with up to 164 different operating power levels. The framework is further applied to engine clusters with different or the same type of engines. We investigated clusters with 2, 3 or 4 different types of engines. The results indicate that the proposed smoothing method is capable of handling discrete set of power modes as high as 164. The results also indicate that the majority of the thrusting occurs during perihelion passages of the intermediate quasi-elliptical orbits, where the power from solar arrays is maximum. However, based on the available power, multiple modes of operation are engaged that may consist of different settings for each engine. The complex engine operation settings is revealed autonomously through the CSC framework.